\documentclass[journal]{IEEEtran}

\usepackage{cite}
\usepackage{todonotes}
\usepackage{graphicx}
\usepackage{subfigure}
\usepackage{epstopdf}
\usepackage{url}
\usepackage{stfloats}
\usepackage{latexsym}
\usepackage{amsfonts}
\usepackage{amsmath}
\usepackage{amssymb}
\usepackage{array}
\usepackage{textcomp}
\usepackage{pxfonts}
\usepackage{txfonts}
\usepackage{cases}
\usepackage{psfrag}
\usepackage{epsfig}
\usepackage{setspace}
\usepackage{enumitem}
\usepackage{xcolor}
\usepackage{multirow}
\usepackage{multicol}
\usepackage{lmodern}

\usepackage[T1]{fontenc}
\usepackage{color}
\usepackage{float}
\usepackage{esint}

\newtheorem{Definition}{Definition}
\newtheorem{Remark}{Remark}
\newtheorem{Lemma}{Lemma}

\newtheorem{Theorem}{Theorem}

\usepackage{todonotes}   

\newcommand{\R}{\ensuremath{\mathbb{R}}}

\newcommand{\bfx}{\mathbf{x}}

\newcommand{\bftheta}{{\pmb{\theta}}}

\begin{document}
\title{\LARGE Identifiability of generalised Randles circuit models}

\author{S.M.Mahdi Alavi, Adam Mahdi, Stephen J. Payne and David A. Howey
\thanks{S.M.M. Alavi was with the Energy and Power Group, Department of Engineering Science, University of Oxford. He is now with the Brain Stimulation Engineering Laboratory, Duke University, Durham, NC 27710, USA. Email: mahdi.alavi@duke.edu}
\thanks{A. Mahdi and S.J. Payne are with the Institute of Biomedical Engineering, Department of Engineering Science, University of Oxford, Old Road Campus Research Building, Oxford, OX3 7DQ, United Kingdom. Emails:  \{adam.mahdi, stephen.payne\}@eng.ox.ac.uk}
\thanks{D.A. Howey is with the Energy and Power Group, Department of Engineering Science, University of Oxford, Parks Road, Oxford, OX1 3PJ, United Kingdom. Email: david.howey@eng.ox.ac.uk}}
\maketitle

\begin{abstract}
The Randles circuit (including a parallel resistor and capacitor in series with another resistor) and its generalised topology have widely been employed in 
electrochemical energy storage systems such as batteries, fuel cells and supercapacitors, also in 
biomedical engineering, for example, to model the electrode-tissue interface in electroencephalography and baroreceptor dynamics. This paper studies identifiability of generalised Randles circuit models, that  is,  whether the model  parameters can be estimated uniquely from the input-output data. It is shown that generalised Randles circuit models are structurally locally identifiable. The condition that makes the model structure globally identifiable is then discussed. Finally, the estimation accuracy is evaluated through extensive simulations.
\end{abstract}

\begin{IEEEkeywords}
Randles circuit, Identifiability, System identification, Parameter estimation.
\end{IEEEkeywords}



\section{Introduction}\label{sec:introduction}
\begin{figure*}
\centering
\includegraphics[scale=.8]{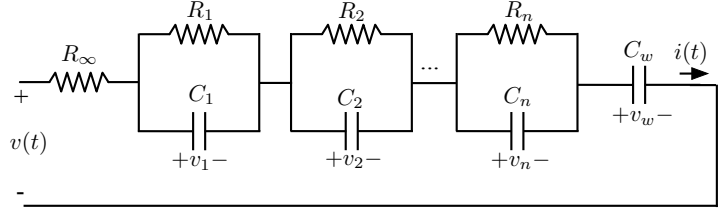}
\caption{The generalised Randles equivalent circuit model.}
\label{Fig:Randles}
\end{figure*}

Randles proposed an equivalent circuit for the kinetics of rapid electrode reactions in \cite{Randles1974}. Since then, the model has been developed and has become the basis for the study of many electrochemical energy storage systems such as batteries, fuel cells and supercapacitors, \cite{Rahn2013, Alavi2015}. Figure~\ref{Fig:Randles} shows a generalised Randles model consisting of an ohmic resistor, $R_\infty$, in series with a number of parallel resistors and capacitors, and a capacitor $C_w$. In electrochemical applications, the ohmic resistor $R_\infty$ represents usually the conduction of charge carriers through electrolyte and metallic conductors. The resistors and capacitors in the parallel pairs represent the charge transfer resistance and the double layer capacitance, respectively, or are an approximation of a diffusion process, \cite{Barsoukov2005, Buller2005}. The number of parallel $R$'s and $C$'s depends on how many of these pairs are required such that the frequency response of the generalised Randles model fits with the device impedance spectra within the frequency range of interests \cite{Andre2011, Hu2011, Birkl2013}. For instance, in \cite{Birkl2013}, the number of the parallel pairs was determined by minimising the error between the model and measured voltages. The capacitor $C_w$, also known as the Warburg term, accounts for a diffusion process, \cite{Barsoukov2005, Buller2005}; or in a battery or supercapacitor, it may represent state of charge \cite{Rahimi-Eichi2014}. It should be noted that the open-circuit voltage source is not considered in the generalised Randles model Figure~\ref{Fig:Randles} and it is only focused on the impedance model in this paper.

The generalised Randles circuit models have also been employed in biomedical engineering. An electric circuit model of the electrode-tissue interface in electroencephalography includes two parallel R-C pairs in series with two resistors and two voltage sources \cite{Mihajlovic2015}, which is a special case of the circuit given in Figure~\ref{Fig:Randles}. 
It is noted that $C_w$ can be considered as a parallel $R$ and $C$ with $R=\infty$. The viscoelastic analog of the generalised Randles model is employed in cardiovascular and cerebral haemodynamics modeling, in describing the viscoelastic properties of the aortic wall and the coupling of the nerve endings of the baroreceptor neurons in the carotid sinus and aortic arch \cite{BugCowBea10, MahStuOttOlu13}; and in relating the fluctuation of the arterial blood pressure with the cerebral blood flow velocity \cite{Mader2014}.

The identification or parameter estimation of the generalised Randles model (with different numbers of parallel R's and C's) is important in condition monitoring, fault diagnosis and control  \cite{Alavi2015, Lee2006, Moubayed2008, Jang2011, Jiang2011, Pattipati2011, Rahmoun2012, Al-Nazer2012, Rahimi-Eichi2014}. In \cite{Lee2006}, the authors showed that the Randles circuit can be used for monitoring the battery charge transfer overvoltage. In \cite{ Moubayed2008}, identification tests for parameter estimation of lead-acid batteries are suggested.  In \cite{Al-Nazer2012}, identification through fitting the impedance spectra in the frequency domain is presented. 

The objective of this paper is to study the identifiability of the generalised Randles model shown in Figure~\ref{Fig:Randles}, that is, whether the model parameters can be estimated uniquely from input-output data. Typically the identifiability problem is divided into two broad areas: parameter estimation accuracy and structural identifiability. \emph{Parameter estimation accuracy} considers the practical aspects of the problem that come with real data such as noise and bias \cite{Raue2009}. In studying \emph{structural identifiability}, on the other hand, one assumes that noise-free informative data is available and therefore it is, in fact, a data-independent concept. Unidentifiable parameters can be assigned an infinite number of values yet still lead to identical input-output data. Thus, structural identifiability is a necessary condition for identifiability and parameter estimation. A number of analytical approaches to structural identifiability have been proposed, including Laplace transform (transfer function) \cite{Cobelli1980, Bellman1970}, Taylor series expansion \cite{Pohjanpalo1978, Chappell1990}, similarity transformations \cite{Vajda1989, Anstett2008, Meshkat2014, Mahdi2014, Glover1974, Distefano1977, VanDenHof1998, Ljung1987}, and differential algebra \cite{Ljung1994,Audoly2001}. In linear systems, it has been shown that controllability and observability properties are closely related to the concept of structural identifiability, \cite{Glover1974,  Distefano1977, VanDenHof1998, Ljung1987}. For example, it was shown that single-input single-output linear time-invariant systems are structurally identifiable if and only if their observer canonical form is controllable (see Chapter 4 in \cite{Ljung1987}). However, controllable and observable systems can still be unidentifiable in the general case \cite{Distefano1977}.

Recently, there has been a significant interest in the identifiability analysis of battery models \cite{Sitterly2011, Rausch2013, Schmidt2010, Forman2012, Moura2013, DAmato2012, Sharma2014, Rothenberger2014}. In \cite{Sitterly2011}, the structural identifiability of a five-element equivalent circuit model (ECM) including two capacitors was discussed by comparing the number of unknown parameters of the Transfer Function (TF) and the circuit. In \cite{Rausch2013}, the structural identifiability of a more general nonlinear ECM is analysed based on the observability conditions. It is shown that cells with serial connections are not observable, demonstrating that complete estimation of the state-of-charge and model parameters from lumped measurements with series cells is not possible trough independent measurements. However, it is shown that  lumped models with parallel connectivities are observable provided that none of the parallel cells are identical.
The identifiability of battery electrochemical models was discussed in \cite{Schmidt2010, Forman2012, Moura2013, DAmato2012}. In particular, it was shown that some of the electrochemical model parameters are not identifiable given typical charge-discharge cycles \cite{Schmidt2010, Forman2012, Moura2013}. In \cite{DAmato2012}, it was shown that the shape of the charge-discharge cycles plays a crucial role in the identifiability of battery parameters. It was also demonstrated that the system identification method can be employed in monitoring the battery film growth. The identifiability problem of Randles ECMs has been studied in \cite{Sharma2014, Rothenberger2014}, where the authors consider models that include up to two capacitors, and the analysis is based on the Fisher information matrix (FIM). The FIM provides some information about the sensitivity of the measurement to the model parameters by using likelihood functions, \cite{Peeters1999}. In \cite{Sharma2014}, a bound of estimation errors was developed by using the Cram\'{e}r-Rao theorem. In \cite{Rothenberger2014}, a method was proposed to optimally shape the battery cycles and improve the identifiability.

This paper shows that the generalised Randles model in Figure~\ref{Fig:Randles} is structurally globally identifiable for $n=1$, and structurally locally identifiable for any finite $n>1$, and becomes globally identifiable assuming an ordering through the generalised circuit. Finally, the identifiability of the model is assessed through extensive simulations.





\section{The model parameterisation and problem statement}
\label{sec: modelling and problem statement}
The state-space and  TF parameterised models of the generalised Randles circuit given in Figure~\ref{Fig:Randles} are derived.

Define the electric current as the system input $u(t)=i(t) \in \R$, the terminal voltage as the system output $y(t)=v(t)\in \R$, the voltages across the internal R-C pairs as the states $\bfx(t)=\left[\begin{array}{cccc} v_1(t) & \cdots & v_n(t) & v_w(t) \end{array}\right]^\top \in \R^{n+1}$, and the model parameters as
\begin{align}
\label{ThetaVec}
\bftheta=\left[\begin{array}{cccccccc}
R_\infty& R_1 & \cdots & R_{n} & C_1 & \cdots & C_n & C_w \end{array}\right],
\end{align}
where $\bftheta$ belongs to some open subset $\mathcal{D}\subset\R^{2n+2}$.

\subsection{State-space parametrisation}
By using Kirchhoff's laws, a state-space model structure of the Randles circuit, parameterised by $\bftheta$, can be written as:
\begin{eqnarray}\label{SS:c}
\left\{\begin{aligned}
&\frac{d}{dt}\bfx(t)  = A(\bftheta)\, \bfx(t) + B(\bftheta)\, u(t)\\
&y(t)    = C (\bftheta)\, \bfx(t) + D(\bftheta)u(t)
\end{aligned}\right.
\end{eqnarray}
where 
$A(\bftheta)\in \R^{(n+1)\times (n+1)}$, $B(\bftheta)\in \R^{n+1}$, $C(\bftheta)^{\top}\in\R^{n+1}$ and $D(\bftheta)\in\R$ are matrices that depend on the parameter vector $\bftheta$, and are given by

\begin{align}
\nonumber
& A(\bftheta) = \left[\begin{array}{ccccc}
-a_1(\bftheta)& 0 & \cdots & 0 & 0 \\
0 & -a_2(\bftheta) & \cdots & 0 & 0\\
\vdots & \vdots & \ddots & \vdots & \vdots\\
0 & 0 & \cdots & -a_{n}(\bftheta) & 0\\
0 & 0 & \cdots &   0 & 0
\end{array}\right],\\ \label{ABCD}
& B(\bftheta) =
\left[\begin{array}{c}
b_1(\bftheta)  \\
b_2(\bftheta) \\
\vdots \\
b_{n}(\bftheta)\\
b_w(\bftheta)\\
\end{array}\right],~
C(\bftheta)^{\top} =
\left[\begin{array}{c}
1  \\
1 \\
\vdots \\
1\\
1\\
\end{array}\right],~ D(\bftheta) = d,
\end{align}
with time constants $\tau_i=R_iC_i$ for $i=1,\cdots,n$ and
\begin{align}
\label{Param:c}
a_i=\frac{1}{\tau_i},~b_i=\frac{1}{C_i},~
b_w=\frac{1}{C_w},~ d=R_\infty.
\end{align}
\subsection{Transfer function parametrisation}
Denote the model's TF by $T(s,\bftheta)$, where $s$ represents the Laplace operator. Using the formula $T(s,\bftheta)=C(sI-A)^{-1}B+D$, a parameterised TF of the generalised Randles circuit can be written as:  \begin{align}\label{TF:gen}
T(s,\bftheta)=\displaystyle \sum_{i=1}^{n}\frac{b_i(\bftheta)}{s+a_i(\bftheta)}+ \frac{b_w(\bftheta)}{s}+D(\bftheta).
\end{align}

\subsection{Problem statement}
Determine the conditions for the model structure \eqref{TF:gen} (equivalently in the state-space form \eqref{SS:c} and \eqref{ABCD}), where  $\bftheta$ is an  unknown parameter vector \eqref{ThetaVec}, to be locally and/or globally structurally identifiable.

\section{Main results}
\label{sec: Main Results}

Following \cite{Ljung1987} a definition of structural identifiability is given as follows.

\begin{Definition}\label{def:ide}
Let  $\mathcal{M}$ be a  model structure with the TF $T(s,\bftheta)$, parametrized by $\bftheta$, where $\bftheta$ belongs to an open subset $\mathcal{D}_T\subset \R^m$, and  consider the equation
\begin{equation}\label{TFd:ide}
T(s,\bftheta)=T(s,\bftheta^\ast)\qquad \text{for almost all } s,
\end{equation}
where $\bftheta, \bftheta^\ast\in\mathcal{D}_T$. Then,  the model structure $\mathcal{M}$ is said to be
\begin{itemize}
\item[-] \emph{globally identifiable} if \eqref{TFd:ide} has a unique solution in $\mathcal{D}_T$,
\item[-] \emph{locally identifiable} if \eqref{TFd:ide} has a finite number of solutions in $\mathcal{D}_T$,
\item[-] \emph{unidentifiabile} if \eqref{TFd:ide} has a infinite number of solutions in $\mathcal{D}_T$.
\end{itemize}
\end{Definition}

\begin{Remark}\upshape\label{Rem:cm}
Instead of Definition~\ref{def:ide}, one can use the so-called coefficient map defined as follows, \cite{Meshkat2014, Mahdi2014}. Consider the monic\footnote{In a monic TF, all coefficients are normalised such that the coefficient of the greatest order in the denominator is 1.} TF:
\begin{equation}\label{eq:Trem}
T(s,\bftheta) = \frac{c_0 (\bftheta)+ c_1(\bftheta)s+ \dots+c_{k_1}(\bftheta)s^{{k_1}}}{d_0 (\bftheta) + d_1(\bftheta)s+ \dots+d_{k_2-1}(\bftheta)s^{{k_2-1}}+s^{k_2}},
\end{equation}
and associate with it the following \emph{coefficient map}  $\mathcal{C}_T:\R^m \supset\mathcal{D}_T \to \R^{k_1+k_2+1}$ defined as
\begin{equation}\label{eq:cm}
\mathcal{C}_T(\bftheta) = \Big[ c_0(\bftheta),\ldots, c_{k_1}(\bftheta),d_0(\bftheta),\ldots, d_{k_2-1}(\bftheta)\Big].
\end{equation}
The model structure $\mathcal{M}$ is globally identifiable if the coefficient map $\mathcal{C}_T$ is one-to-one (injective); locally identifiable if $\mathcal{C}_T$ is many-to-one; and  unidentifiable if it is infinitely many-to-one.
\end{Remark}

\smallskip

The following lemma will be used in the proof of the identifiability of the generalised Randles circuit model.

\begin{Lemma}\label{lemma}
Let $C_T(m)$ be a coefficient map associated with the TF 
\[
T(s,\bftheta)=\displaystyle \sum_{i=1}^{m}\frac{b_i}{s+a_{i}}+d,
\]
where $\bftheta=(a_1,\dots,a_m,b_1,\ldots,b_m,d)$ and $a_i\in\R$ for $i=1,\ldots, m$ are pairwise different. Then the following statements hold:
\begin{itemize}
\item[(a)] If $m=1$, then $C_T(m)$ is one-to-one.
\item[(b)] If $m>1$, then $C_T(m)$ is many-to-one.
\end{itemize}
\end{Lemma}

{\bf Proof.} {\it Part(a).} For $m=1$, the identifiability equation \eqref{TFd:ide} is given by
\[
\frac{b_1}{s+a_1}+d = \frac{b^\ast_1}{s+a^\ast_1}+d^\ast
\]
and has a unique solution $(a_1,b_1,d)=(a^\ast_1,b^\ast_1,d^\ast)$, which proves part (a).

{\it Part(b).} For $m>1$, the identifiability equation \eqref{TFd:ide} is given by
\begin{equation}\label{ie:dn}
\displaystyle \sum_{i=1}^{m}\frac{b_i}{s+a_{i}}+d = \displaystyle \sum_{i=1}^{m}\frac{b^\ast_i}{s+a^\ast_{i}}+d^\ast.
\end{equation}
We claim that equation \eqref{ie:dn} admits only finite (more precisely, $m! = 1\times2\times \cdots \times m$) number of solutions. To prove the claim note that \eqref{ie:dn} is the equality of two rational functions, which is satisfied provided that
\begin{equation}\label{lemma:den}
(s+a_1)\times \cdots \times (s+a_m) = (s+a_1^\ast)\times \cdots \times (s+a_m^\ast).
\end{equation}
Since $m$ distinct roots uniquely characterise a monic polynomial of degree $m$, and there are $m!$ permutations of $n$ roots of $(s+a_1)\times \cdots \times (s+a_m)$, equation \eqref{lemma:den} has $m!$ solutions. Now, let us fix the permutation $s:\{1,\ldots, m\}\to \{1,\ldots,m\}$ and consider an assignment  $(a_1,\ldots,a_m) = (a_{s_1}^\ast,\ldots,a_{s_m}^\ast)$.  Since $a_i$ for $i=1,\ldots,m$ are assumed to be pairwise distinct, the expressions $1/(s+a_i)$, thought of as functions of the variable $s$, are linearly independent. Finally, since each side of equation \eqref{ie:dn} is a linear combination of linearly independent functions, we immediately obtain that $b_j=b_{s_j}^\ast$ for $i=1,\ldots,m$. This concludes the proof of part (b).\hfill{$\square$}



\smallskip

Now, we introduce the concept of reparametrisation (see e.g. \cite{Meshkat2014}), which we will use the proof of our main theorem.
\begin{Definition}\upshape\label{def:Rep}
A \emph{reparametrisation} of the model structure $\mathcal{M}$  with the coefficient map $\mathcal{C}_T$ is a map $\mathcal{R}:\R^k \supset\mathcal{D}_T \to \R^{m}$ such that
\begin{equation}
{\rm Im}\,(\mathcal{C}_T \circ \mathcal{R}) = {\rm Im}\, \mathcal{C}_T,
\end{equation}
where ${\rm Im}$ denotes the image of the map. Moreover, the reparameterisation is identifiable if the map $\mathcal{C}_T \circ \mathcal{R}: \R^k \supset\mathcal{D} \to \R^{k_1+k_2+1}$ is identifiable.
\end{Definition}

\smallskip

The main result of this section is the following theorem which describes the identifiability of the generalised Randles circuit.
\begin{Theorem}\label{Thm:main}
Let $\mathcal{M}_{RC}(n)$ denotes the state-space model structure \eqref{SS:c} with the matrices \eqref{ABCD} parametrised by \eqref{Param:c}, where $n$ is the number of parallel $RC$ elements connected in series (see Figure~\ref{Fig:Randles}). Then the following conditions hold:
\begin{itemize}
\item[(a)] If $n=1$, then the model structure  $\mathcal{M}_{RC}(n)$ is globally identifiable.
\item[(b)] If $n>1$, then the model structure  $\mathcal{M}_{RC}(n)$ is locally identifiable.
\item[(c)] If $n>1$, and there is an ordering through the generalised circuit as
    \begin{align}
    \label{EIS condition1}
    & a_n < a_{n-1} < \cdots < a_1,
        \end{align}
then the model structure  $\mathcal{M}_{RC}(n)$ is globally identifiable.
\end{itemize}
\end{Theorem}

{\bf Proof.} Consider the model structure $\mathcal{M}_{RC}(n)$ and let $T(s,\bftheta)$ denote the corresponding  TF given by \eqref{TF:gen}, where the parameter vector $\bftheta$ is given by \eqref{ThetaVec}. We write $T(s,\bftheta)$ as a rational function \eqref{eq:Trem} of degree $k_1=k_2=n+1$ and associate with it the coefficient map $\mathcal{C}'_T:\R^{2n+2}\supset\mathcal{D}\to \R^{2n+3}$.



{\it Part (a)}. For $n=1$, the coefficient map $\mathcal{C}'_T:\R^{4}\supset \mathcal{D}\to R^{5}$ can be written explicitly as
\[
\mathcal{C}'_T(\bftheta) = \Big[ c_0(\bftheta), c_1(\bftheta), c_2(\bftheta),d_0(\bftheta),d_1(\bftheta)\Big],
\]
where
\[
\begin{aligned}
&c_0(\bftheta) = \frac{1}{R_1C_1C_w},\quad   c_1(\bftheta) = \frac{R_{\infty}C_w+R_1C_w+R_1C_1}{R_1C_1C_w}, \\ & c_2(\bftheta) = R_{\infty},\quad 
d_0(\bftheta) = 0,\quad   d_1(\bftheta)=\frac{1}{R_1C_1}.
\end{aligned}
\]
By direct computation, we can check that equation $\mathcal{C}'_T(\bftheta)= \mathcal{C}'_T(\bftheta^*)$  admits a unique solution $\bftheta=\bftheta^*$. Thus the coefficient map is one-to-one, and  the model structure is globally identifiable. 

{\it Part (b).} For $n>1$, the coefficient map $\mathcal{C}'_{T}$ can be written as the following composition
\begin{equation}\label{CT:com}
\mathcal{C}'_{T} = \mathcal{C}_T \circ \mathcal{R}^c,
\end{equation}
where the map $\mathcal{R}^c:\R^{2n+2}\supset \mathcal{D}\to R^{2n+2}$ is the reparametrisation $\mathcal{R}^c(\bftheta) = \big(a_1,\ldots,a_n,b_1,\ldots,b_n,b_w,d\big)$ defined by \eqref{Param:c}; and $\mathcal{C}_T$ is the coefficient map associated with the TF
\begin{equation}\label{eq:TF_ct}
T(s, \bftheta_{a,b}) = \sum_{i=1}^n \frac{b_i}{s+a_i}+\frac{b_w}{s}+d,
\end{equation}
where $\bftheta_{a,b} = (a_1,\ldots,a_n,b_1,\ldots,b_n,b_w,d)\in\R^{2n+2}$ . By Lemma~\ref{lemma} (for $m=n+1$ and $a_{n+1}=0$) the map $C_T$ is many-to-one; and $\mathcal{R}^c$ is a one-to-one map with an inverse
\[
R_{\infty} = d,\,\, C_w=\frac{1}{b_w},\,\, R_i = \frac{b_i}{a_i}, \,\, C_i = \frac{1}{b_i}, \,\,i=1,\ldots,n.
\]
The map  $\mathcal{C}'_T$ is many-to-one, since it is a composition of a one-to-one with many-to-one map. Thus for $n>1$, the model structure $\mathcal{M}_{RC}(n)$ is locally identifiable.

{\it Part (c).} Finally, the identifiability equation \eqref{ie:dn} under the condition \eqref{EIS condition1} admits a unique solution, which concludes part (c).\hfill{$\square$}



\begin{Remark}\upshape
The same procedure is applicable to the discrete time model. The Euler's first order approximation is the simplest approximation.  Its identifiability analysis is easier because the coefficients of the discrete time TF equals the number of parameters. However, it might lead to numerical instability. If higher order approximations are applied, the challenge remains to prove whether there is a finite-to-one maps between the discrete time model and parameters.  
\end{Remark}

\section{Simulations}
\label{sec: Simulations}

Consider a 6-element Randles circuit with $R_{\infty}$, $C_w$ $R_i$, $C_i$, $i=1,2$. In this section, accuracy of the estimation in the presence of noise-free and noisy data, subject to random initial guess of estimations is studied. First the informative data set used in the simulations and the way it is generated are described. More details are provided in references \cite{Ljung1987, Soderstrom1989, Norton1986, Zhu2001}.

\begin{figure*}
\centering
\subfigure[Excitation input $u(t)$.]{\label{fig:u}
\includegraphics[scale=.45]{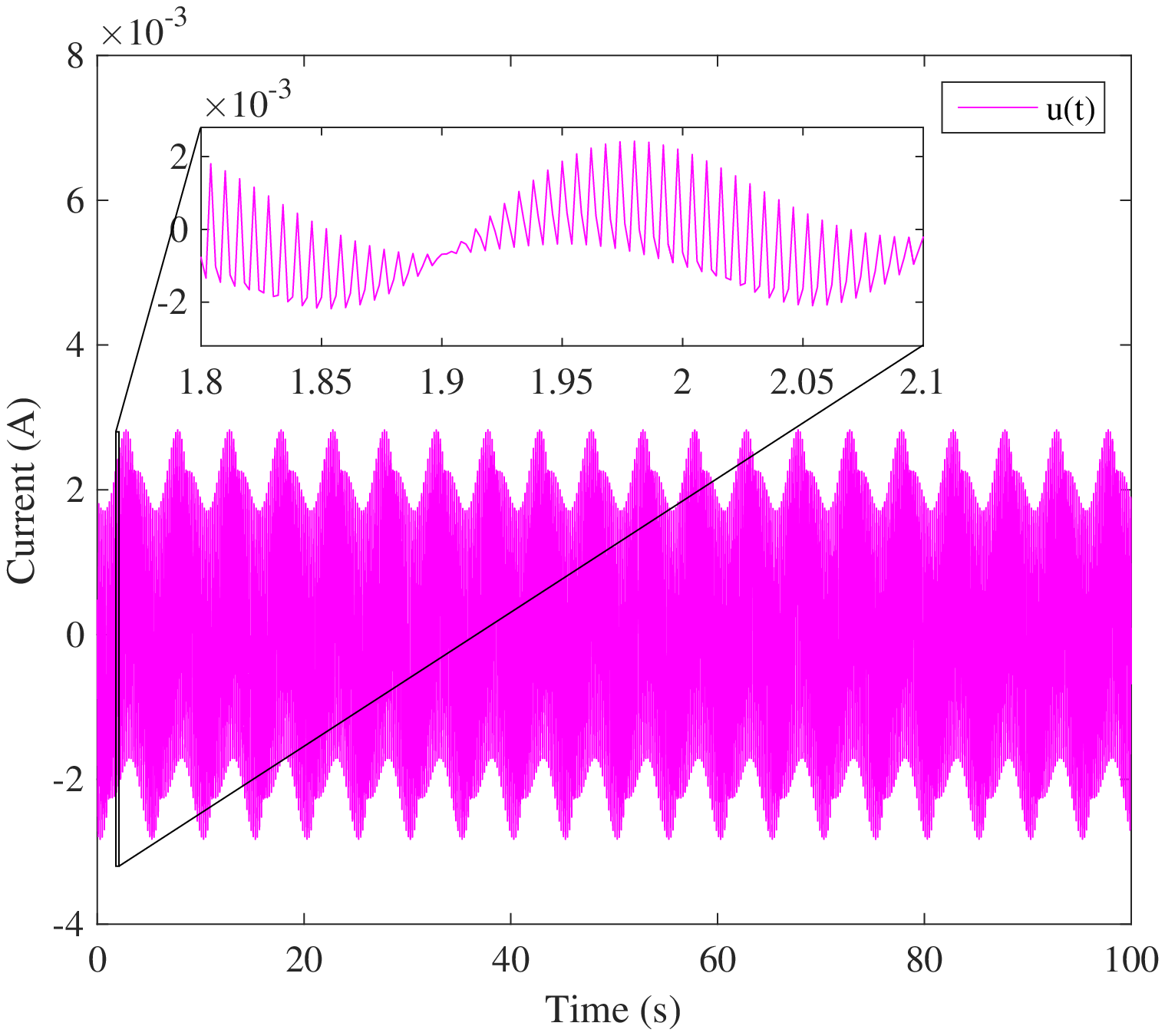}}
\subfigure[Associated output $y(t)$.]{\label{fig:y}
\includegraphics[scale=.45]{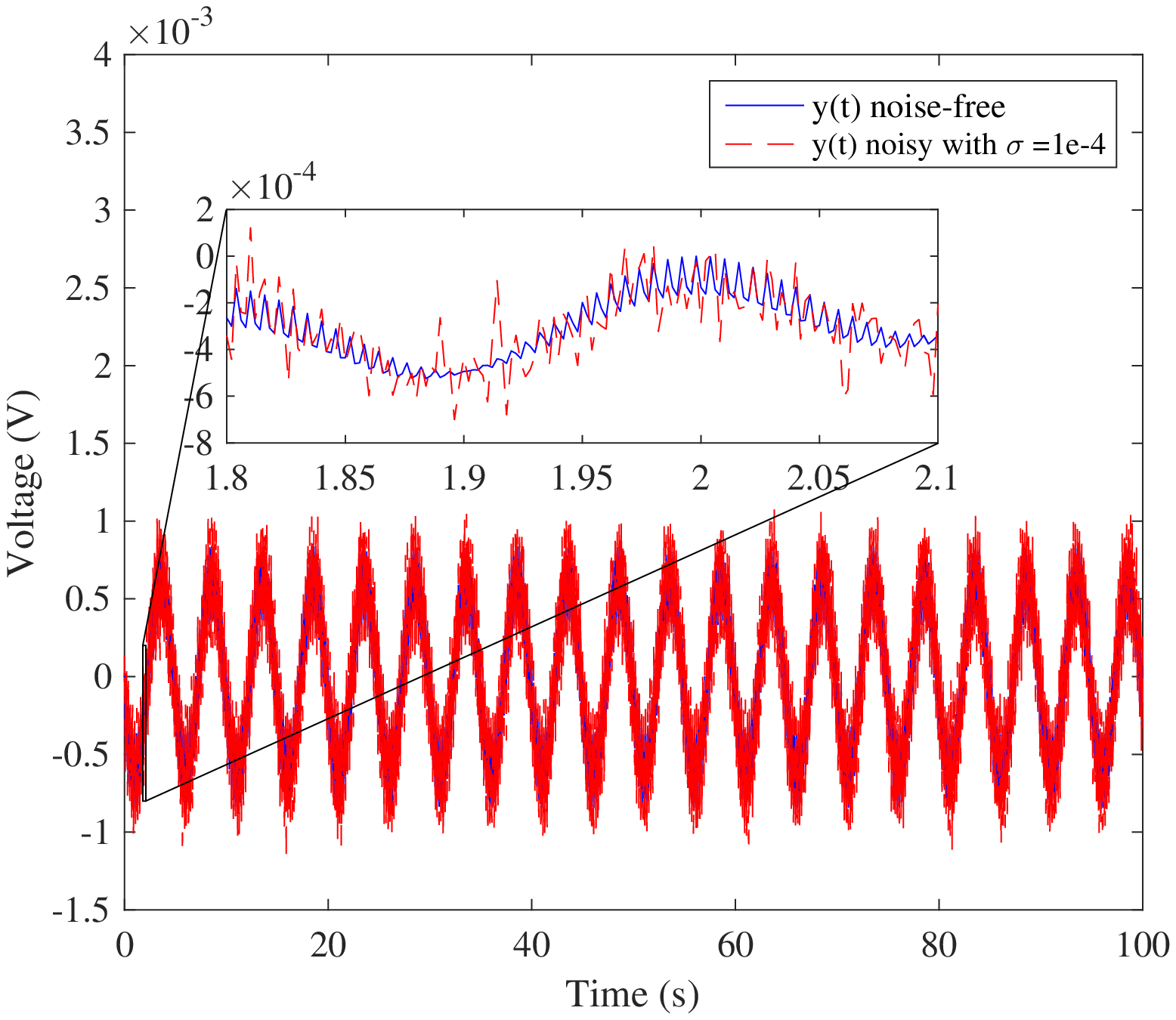}}
\caption{The multi-sine input excitation signal and its associated output for $R_\infty=0.05~\Omega,~~ R_1=0.2~\Omega,~~ C_1=0.3~\mbox{F},~~ R_2=0.4~\Omega,~~ C_2=0.6~\mbox{F},$ and $C_w=300~\mbox{F}$. The voltage response, Figure \ref{fig:y}, is computed using the model \eqref{SS:c} in the noise-free and noisy case with standard deviation $\sigma=10^{-4}$. The DC offset has been removed to have zero-mean signals. }
\label{uydata}
\end{figure*}

\subsection{Generation of informative data set}
A data set is informative if the input is persistently exciting. A persistently exciting input adequately excites all modes of the system. In linear systems, the order of the system determines the order of persistent excitation. The order of a persistent excitation equals the number of coefficients of the monic TF that need to be identified, (see Theorem 13.1. in \cite{Ljung1987}). 
The monic TF of the 6-element Randles circuit is
\begin{align}
\label{ex-tf}
H(s)=\frac{f_{3}s^3+f_{2}s^2+f_{1}s+f_{0}}{s^3+g_{2}s^2+g_{1}s+g_{0}}
\end{align}
with unknown coefficients to be identified. Therefore, the necessary order of the persistently exciting input is 7. In the frequency domain, this means that the spectrum of the excitation input should have at least 7 nonzero points. This simulation focuses on multi-sine excitation signals which are widely employed in electrochemical impedance spectroscopy (EIS) techniques, \cite{Alavi2015}, \cite{Barsoukov2005}, \cite{Howey2014a}. The same procedure can be applied to other inputs. The multi-sine signal is given by
\begin{align}
\label{multisins}
u(t)=\displaystyle \sum_{j=1}^l m_j ~\cos(\omega_jt+\phi_j),
\end{align}
where $l$ represents the number of sinusoids, and $m_j$, $\omega_j$ and $\phi_j \in [-\pi, \pi)$ denote the magnitude, frequency in radians per second and phase in radians, respectively. The spectrum of the multi-sine signal is given by
\begin{align}
\label{multisins-spec}
\Phi_u(\omega)=2\pi\displaystyle \sum_{j=1}^l \frac{m_j^2}{4} [\delta(\omega-\omega_j)+\delta(\omega+\omega_j)],
\end{align}
where $\delta(\omega)$ is the delta function or impulse at frequency $ \omega_j$. The spectrum of each sinusoid signal contains two nonzero points, therefore 4 sinusoid signals are enough to generate the informative data set for the 6-element Randles circuit model.

The magnitudes $m_j$, frequencies $\omega_j$ and phases $\phi_j$ are arbitrary real values. Specific applications might impose additional constraints. 
For instance, in EIS techniques, the magnitude of the input signal may vary from milliampere to ampere depending on the size of the energy storage system. For the sake of simplicity, all magnitudes are assumed equal to $m_j=10^{-3}$ for $j=1,\cdots,4$. The frequencies could be equally or logarithmically spread over the frequency band $\omega_{\min} \leq \omega \leq \omega_{\max}$. The values of $\omega_{\min}$ and $\omega_{\max}$ depend on the system dynamics. If $m_j$'s are equal, the Schroeder phase choice is suggested to reduce the Crest factor \cite{Schroeder1970, Ljung1987}.

\begin{table*}[h]
\caption{\normalsize Estimation results.}
\centering
\begin{tabular}{|l|l|l|l|l|l|l|l|}
\hline
\multirow{2}{*}{} & Parameter& $R_\infty$  & $R_1$  & $C_1$ & $R_2$  & $C_2$  & $C_w$ \\
                 &   True Value& 0.05          & 0.2       &    0.3     &    0.4     &    0.6    & 300\\ \hline
\multirow{3}{*}{Noise free (9 outlier estimates)} & mean&0.0552   & 0.1883 &   0.2919    &0.4061 &   0.5866 & 289.3709
\\
                 &   st.d. & 0.0010   & 0.0393  &  0.0617  &  0.0360  &  0.0734 &  48.3550
 \\ 
                  &  $e_r$&     10.36 &   5.87 &   2.68   & 1.52  &  2.22 &   3.54 \\ \hline
\multirow{3}{*}{Noise $\sigma=10^{-4}$ (11 outlier estimates)}& mean  & 0.0552& 0.1847  &  0.2886 &   0.4094  &  0.5803 & 300.9287\\ 
                  & st.d. & 0.0012  &  0.0451   & 0.0706 &   0.0415   & 0.0867  & 87.8931\\ 
                  & $e_r$ &  10.38  &  7.63  &  3.79  &  2.34  &  3.28  &  0.31 \\ \hline
\end{tabular}
\label{Table-results}
\end{table*}


\begin{Remark}
A multi-sine signal may not exactly be zero mean. Sometimes a  uni-directional current flow is required. In this case, the multi-sine excitation current needs to be superimposed on a known constant DC offset current, (see examples in \cite{Alavi2015} and \cite{Howey2014a}). In system identification, the data are typically pre-filtered to remove these types of offsets.
\end{Remark}


\subsection{The results discussion}
The following operating point is arbitrarily selected \cite{Alavi2015, Barsoukov2005}:
\begin{align*}
& R_\infty=0.05~\Omega,~~ R_1=0.2~\Omega,~~ C_1=0.3~\mbox{F}\\ & R_2=0.4~\Omega,~~ C_2=0.6~\mbox{F},~~ C_w=300~\mbox{F}
\end{align*}
The multi-sine excitation signal, Figure \ref{fig:u}, is generated by using model \eqref{multisins} with $l=4$, $m_j=10^{-3}$, $f_{\min}=0.2$ Hz, $f_{\max}=500$ Hz and Schroeder phase with $\phi_1$ randomly chosen at $1.9775$.  The Crest factor of the signal is  1.9999. The voltage response, Figure \ref{fig:y}, has been computed using the model \eqref{SS:c}. The DC offset has been removed to ensure zero-mean signals.

Given the true values, the smallest time constant is $\tau_{\min}=R_1C_1=0.06$ s. The sampling time should be several times larger than $1/\tau_{\min}$. In this simulation, the sampling frequency is chosen at  $f_s=500$ Hz, which is 16.6 times greater than the inverse of the minimum time constant. The test duration should typically be several times larger than the maximum time constant, which is $\tau_{\max}=R_{\infty}C_w=15$ s. In \cite{Zhu2001}, 6 to 8 times $\tau_{\max}$ is suggested, however, this might vary in different applications. In this simulation, a test duration of $t_d=100$ s is applied. The continuous-time TF \eqref{ex-tf} is identified using Matlab's system identification toolbox \cite{Ljung1988}. The circuit parameters are calculated directly from the coefficients of the TF using the formulas shown in Table II and discussed later. In order to study the consistency of results, each test is repeated 100 times, every run with a random initial guess of parameters. The roots of the denominator must be positive real numbers. Those estimations, which lead to complex or negative poles for the TF, or to $C_w > 1000$ or to $C_i>10$, $i=1,2$, are considered outliers and discarded from the analysis. 

The relative mean error, $e_{r}$ is defined as follows:
\begin{align}
\label{mean-error}
e_{r}=100\times \left|\frac{\mbox{true~value}-\mbox{mean~of~estimations}}{\mbox{true~value}}\right|\%.
\end{align}

\begin{figure}[htp]
\centering
\includegraphics[scale=.45]{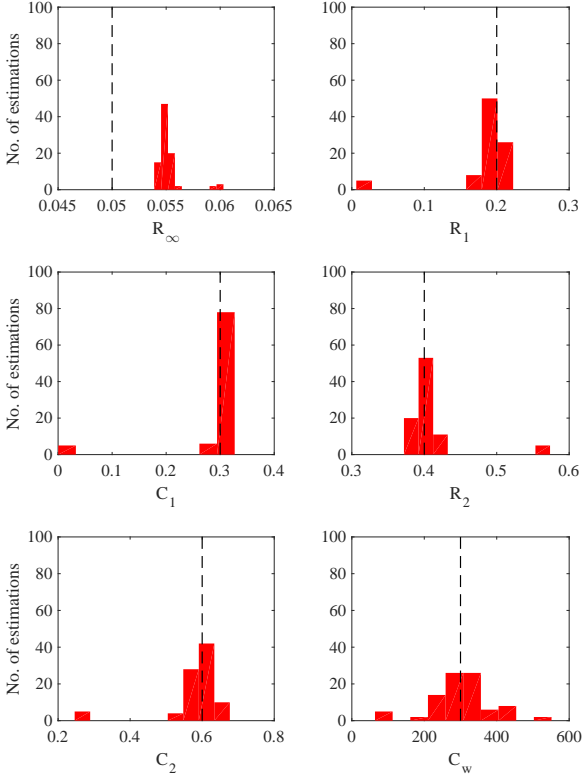}
\caption{Histogram of accepted estimations in 100 runs from noisy data.}
\label{hdata}
\end{figure}

The estimation accuracies from noise-free and noisy data with a zero-mean noise and standard deviation $\sigma=10^{-4}$ are compared together. The signal to noise ratio is $\frac{10^{-3}}{10^{-4}}= 10$. Table \ref{Table-results} shows the mean, standard deviation and relative errors of the estimations. Figure \ref{hdata} shows the histograms of the estimations from the noisy data. 

Regardless of the random initial condition, in more than 80\% of simulations, the noise-free estimations converge with $e_r$ less than 10\% for $R_i$, $C_i$, $i=1,2$ and for $C_w$. The largest $e_r$ is $10.36\%$ for $R_{\infty}$ (Table \ref{Table-results}). 

Figure \ref{hdata} shows that estimations from noisy data are more distributed around the true values. However as it is seen in Table \ref{Table-results}, the mean values of estimations in both noise-free and noisy cases remain the same. The relative mean errors of estimations are again less than 10\% for all parameters except of $R_{\infty}$, which is again the largest and $10.38\%$. Also, $e_r$'s increase with noisy data compared to that of in the noise-free case.  

The largest standard deviation is for $C_w$. This could be because of the pure integrator associated with $C_w$, which appears in the TF, which might require some modifications of the data set, \cite{Anderson2005}. For instance, a methodology has been proposed in \cite{Sitterly2011} that removes the integral term by modifying the input signal as $i_{\mbox{\scriptsize modified}}(t)=\int i(t) dt$. 

\subsection{Calculating the coefficients of $R-R||C-R||C-C$ circuit}
Here we show how to compute the parameters of the 6-element Randles circuit. Using (\ref{TF:gen}), the circuit's TF is:
\begin{align*}
H(s)=\frac{f_{3}s^3+f_{2}s^2+f_{1}s+f_{0}}{s^3+g_{2}s^2+g_{1}s+g_{0}}
\end{align*}
The relationships between the coefficients and the circuit parameters are given by:
\begin{eqnarray}\label{Param:invmap}
\begin{aligned}
& f_{3}=d,~ f_{2}=b_1+b_2+b_w +(a_1+a_2)d\\
& f_{1}=a_2 b_1+a_1 b_2+(a_1+a_2)b_w +a_1 a_2 d,~ f_{0}=a_1a_2b_w\\
& g_{2}=a_1+a_2,~g_{1}=a_1a_2,~g_{0}=0
\end{aligned}
\end{eqnarray}
where the parameters $a_i$, $b_i$, $b_w$ and $d$ are as defined in (\ref{Param:c}).

Because the circuit has an integrator, the identification method should be set up such that a pole of the denominator is fixed at $s=0$. The identification software typically allows to fix a number of poles and zeros at certain values, \cite{Ljung1988}. Using the first equation of (\ref{Param:invmap}) and (\ref{Param:c}), $R_\infty$ is given by:
\begin{align*}
R_\infty =f_{3}
\end{align*}

The roots of $s^3+g_{2}s^2+g_{1}s+{g}_0=0$ are ${a}_1$ and ${a}_2$, and the one that has been fixed at $s=0$. By using the condition \eqref{EIS condition1}, select the smallest root as ${a}_2$, and the remaining root as ${a}_1$. From (\ref{Param:c}), the circuit's time constants are obtained as follows:
\begin{align*}
& \tau_i = \frac{1}{{a}_i} \mbox{~for~} i=1,2.
\end{align*}

From (\ref{Param:invmap}) and \eqref{EIS condition1}, $b_1$, $b_2$ and $b_w$ are obtained by solving the following set of equations for $X$:
\begin{align*}
& \left[ \begin{array}{ccc}
1 & 1 & 1 \\
{a}_2& {a}_1 & {a}_1+{a}_2\\
0 & 0 & {a}_1{a}_2\end{array} \right]X=  \left[ \begin{array}{c}
{f}_2-({a}_1+{a}_2){f}_3\\
{f}_1-{a}_1{a}_2{f}_3\\
{f}_0\end{array} \right]
\end{align*}
where,
\begin{align*}
& X=\left[\begin{array}{ccc}
{b}_1& {b}_2 & {b}_w \end{array}\right]^\top
\end{align*}

Other parameters of the circuit are subsequently obtained as:
\begin{align*}
{C}_w=\frac{1}{{b}_w},\quad {C}_i=\frac{1}{{b}_i}, \quad
{R}_i=\frac{{\tau}_i}{{C}_i},\quad i=1,2.
\end{align*}

The parameters of different topologies can simply be calculated using the same approach. Table \ref{Table:invmap} provides formulas for four widely used Randles models.


\begin{table*}
\caption{\normalsize The coefficients of four widely used Randles circuit models.}
\centering
\begin{tabular}{|c|c|c|c|}
\hline
\includegraphics[scale=.5]{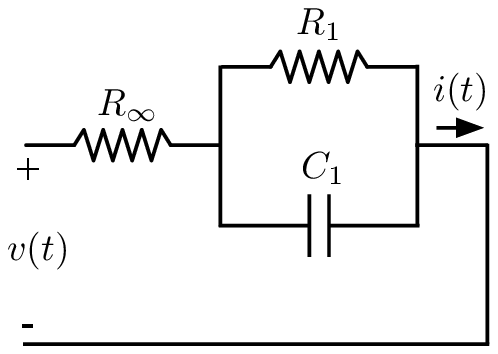} & \includegraphics[scale=.5]{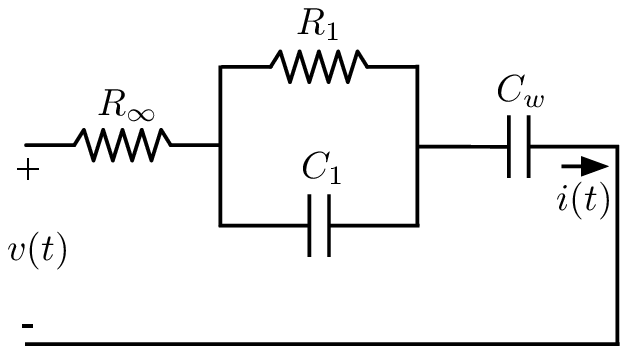} & \includegraphics[scale=.5]{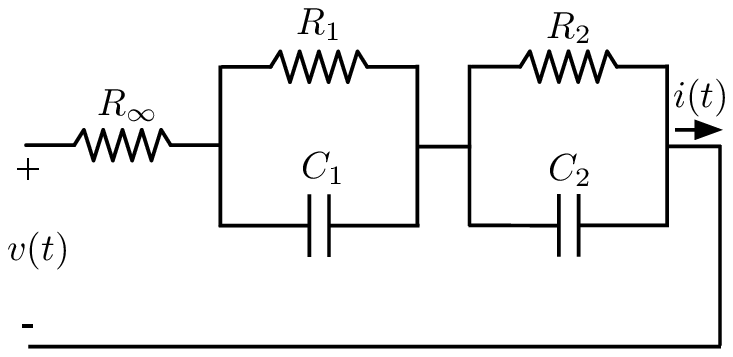} & \includegraphics[scale=.5]{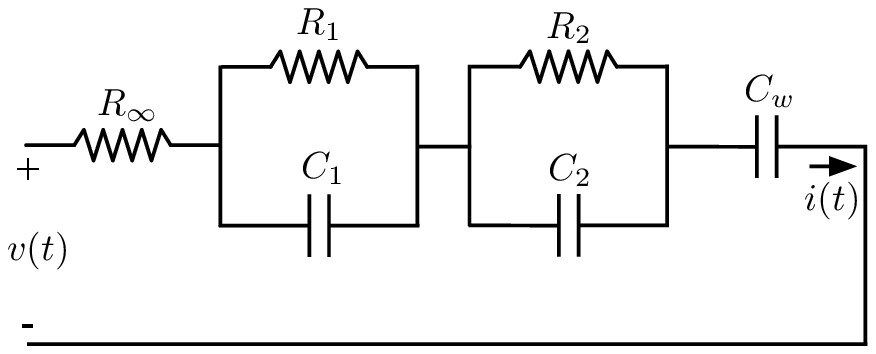}  \\
\hline
&&&\\
\scriptsize Estimated TF: & \scriptsize Estimated TF: & \scriptsize Estimated TF: &\scriptsize Estimated TF:\\
${H}(s)=\frac{{f}_1s+{f}_0}{s+{g}_0}$ & ${H}(s)=\frac{{f}_2s^2+{f}_1s+{f}_0}{s^2+{g}_1s+{g}_0}$ & ${H}(s)=\frac{{f}_2s^2+{f}_1s+{f}_0}{s^2+{g}_1s+{g}_0}$ &
${H}(s)=\frac{{f}_3 s^3 + {f}_2s^2+{f}_1s+{f}_0}{s^3+{g}_2 s^2+{g}_1s+{g}_0}$ \\
& \scriptsize Set up to identify a pole at $s=0$& &  \scriptsize Set up to identify a pole at $s=0$\\
&&&\\
\hline
&&&\\
${R}_\infty ={f}_1$ & ${R}_\infty ={f}_2$ &${R}_\infty ={f}_2$&${R}_\infty ={f}_3$\\
&&&\\
${a}_1={g}_0$&$\mbox{\scriptsize roots}([1~ {g}_1 ~{g}_0]) \Rightarrow {a}_1,0$&$\mbox{\scriptsize roots}([1~ {g}_1 ~{g}_0]) \Rightarrow {a}_1,{a}_2$&$\mbox{\scriptsize roots}([1~ {g}_2~{g}_1 ~{g}_0])\Rightarrow {a}_2,{a}_1,0$\\
& &{\scriptsize choose} ${a}_2 < {a}_1$&{\scriptsize choose} ${a}_2 < {a}_1$\\
&&&\\
${\tau}_1 = \frac{1}{{a}_1} $& ${\tau}_1=\frac{1}{{a}_1} $&${\tau}_i=\frac{1}{{a}_i} $&${\tau}_i=\frac{1}{{a}_i} $\\
&&&\\
&$AX=B$&$AX=B$&$AX=B$ \\
&\scriptsize $A=\left[ \begin{array}{cc}
1 & 1 \\
0 & {a}_1 \end{array} \right]$&\scriptsize $A=\left[ \begin{array}{cc}
1 & 1 \\
{a}_2 & {a}_1 \end{array} \right]$&\scriptsize $A=\left[ \begin{array}{ccc}
1 & 1 & 1 \\
{a}_2& {a}_1 & {a}_1+{a}_2\\
0 & 0 & {a}_1{a}_2\end{array} \right]$ \\
&&&\\

&\scriptsize $B=\left[ \begin{array}{c}
{f}_1-{a}_1{f}_2\\
{f}_0 \end{array} \right]$&\scriptsize $B=\left[ \begin{array}{c}
{f}_1-({a}_1+{a}_2){f}_2\\
{f}_0-{a}_1{a}_2{f}_2 \end{array} \right]$&\scriptsize $B=\left[ \begin{array}{c}
{f}_2-({a}_1+{a}_2){f}_3\\
{f}_1-{a}_1{a}_2{f}_3\\
{f}_0\end{array} \right]$\\
&&&\\
 ${b}_1={f}_0-{a}_1{f}_1$& $X=\left[\begin{array}{cc}
{b}_1&  {b}_w \end{array}\right]^\top$  & $X=\left[\begin{array}{cc}
{b}_1& {b}_2 \end{array}\right]^\top$  & $X=\left[\begin{array}{ccc}
{b}_1& {b}_2 & {b}_w \end{array}\right]^\top $\\
&&&\\
${C}_1 = \frac{1}{{b}_1} $ & ${C}_1= \frac{1}{{b}_1} $ &${C}_i= \frac{1}{{b}_i} $&${C}_i=\frac{1}{{b}_i} $\\
&&& \\
& ${C}_w= \frac{1}{{b}_w} $ &&${C}_w= \frac{1}{{b}_w} $\\
&&& \\
${R}_1 =\frac{{\tau}_1}{{C}_1}$ &${R}_1 =\frac{{\tau}_1}{{C}_1}$&${R}_i =\frac{{\tau}_i}{{C}_i},~i=1,2$&${R}_i =\frac{{\tau}_i}{{C}_i},~i=1,2$\\
&&&\\
\hline
\end{tabular}
\label{Table:invmap}
\end{table*}

\section{Conclusions}
We showed that the generalised Randles circuit model is locally identifiable and the model structure becomes globally identifiable if an ordering through the  circuit is assumed. The results were confirmed through extensive simulations. Finally, explicit formulas for the coefficients of widely used Randles circuits were presented.   

\section*{Acknowledgements}
This work was funded by the University of Oxford EPSRC Impact Acceleration Account Technology Fund Awards EP/K503769/1 and EP/K036157/1. The authors would like to thank the anonymous reviewers and the editor for their fruitful comments that significantly improved the paper. We would like also to extend our thanks to Ross Drummond, Stephen Duncan, Xinfan Lin and Shi Zhao for their feedback on the first version of this work.

\bibliographystyle{IEEEtran}
\bibliography{Ref}

\end{document}